\documentclass[12pt,reqno]{article}

\usepackage{pstricks,pst-node,multido}
\usepackage{amsmath}
\usepackage{amssymb, amscd}

\DeclareMathOperator{\Aut}{Aut}
\DeclareMathOperator{\GL}{GL}

\makeatletter \renewcommand*\l@section{\@dottedtocline{1}{1.5em}{2.3em}} \makeatother

\newcounter{daana} \newcounter{daanb}

\newcommand{\lista}[2]
{\begin{list}{#1}{
\usecounter{daana}
\setlength{\listparindent}{\parindent}
\setlength{\leftmargin}{3.5em}
\setlength{\itemindent}{0mm}
\setlength{\labelsep}{.6em}
\setlength{\labelwidth}{10em}
\setlength{\itemsep}{0em}
\setlength{\parsep}{1mm}
\setlength{\listparindent}{\parindent}
#2}}

\newcommand{\listb}[2]
{\begin{list}{#1}{
\usecounter{daanb}
\setlength{\rightmargin}{\leftmargin}
\setlength{\listparindent}{\parindent}
#2}}

\newtheorem{numb}[equation]{\hspace*{-.35em}}

\newenvironment{prop} {\begin{numb} {\bf Proposition.\ }\sl} {\end{numb}}

\newenvironment{theo} {\begin{numb} {\bf Theorem.\ }\sl} {\end{numb}}
\newenvironment{lemm} {\begin{numb} {\bf Lemma.\ }\sl} {\end{numb}}

\newenvironment{rema} {\begin{numb} {\bf Remark.\ }\rm} {\end{numb}}
\newenvironment{defi} {\begin{numb} {\bf Definition.\ }\rm} {\end{numb}}

\numberwithin{equation}{section}

\newcommand{\daanmath}{\mathbb}
\newcommand{\zz}{\daanmath{Z}}

\newcommand{\rr}{\daanmath{R}}
\newcommand{\cc}{\daanmath{C}}

\newcommand\sur{\mathrel{\to\kern-1.8ex\to}}
\newcommand\iso{\mathrel{\hookrightarrow\kern-1.8ex\to}}

\newcommand\longhookrightarrow{\lhook\joinrel\longrightarrow}

\newcommand\longinj{\longhookrightarrow}
\newcommand\longsur{\mathrel{\longrightarrow\kern-1.8ex\to}}
\newcommand\longiso{\mathrel{\longhookrightarrow\kern-1.8ex\to}}

\newcommand\hugerightarrow{\xrightarrow{\hspace*{2.3em}}}
\newcommand\hugehookrightarrow{\lhook\joinrel\hugerightarrow}
\newcommand\hugeinj{\hugehookrightarrow}
\newcommand\hugesur{\mathrel{\hugerightarrow\kern-1.8ex\to}}
\newcommand\hugeiso{\mathrel{\hugehookrightarrow\kern-1.8ex\to}}

\newcommand{\<}{\langle} 
\renewcommand{\>}{\rangle}
\newcommand{\noi}{\noindent}

\newcommand{\be}{\begin{equation}}
\newcommand{\ee}{\end{equation}}
\newcommand{\kopje}[1]{\bigskip\noi{\bf\large #1.\ \ }\\[.5ex] \noi}

\newcommand{\proof}{\par\noi{\sl Proof.\ \ }}
\newcommand{\BOX}{{}\raisebox{-.05ex}{\makebox[1em][r]{$\Box$}}\bigskip}
\newcommand{\blok}{\hspace*{0cm}\hfill\BOX\par}
\newcommand{\blokk}{\hspace*{0cm}\hfill\BOX\vspace{-.6\baselineskip}\par}
\newcommand{\col}{\text{\upshape :\ }}

\newcommand{\mytitle}{Horizontal configurations of points in link complements}
\title{\mytitle}

\begin{document} 
%


\setlength{\headsep}{8mm}

\begin{center}{\parbox{.77\textwidth}{\center{\huge \mytitle\par}}} \\[5mm] Daan Krammer \\[1mm] 26 October 2004 \end{center}

\thispagestyle{empty}

\begin{abstract} For any tangle $T$ (up to isotopy) and integer $k\geq 1$ we construct a group $F(T)$ (up to isomorphism). It is the fundamental group of the configuration space of $k$ points in a horizontal plane avoiding the tangle, provided the tangle is in what we call Heegaard position. This is analogous to the first half of Lawrence's homology construction of braid group representations. We briefly discuss the second half: homology groups of $F(T)$.
\end{abstract}


\section{Introduction}

In her thesis \cite{law} Ruth Lawrence introduced and studied certain representations of braid groups. She related her representations to the Jones polynomial (see also \cite{big3}). Some of her representations were later shown to be faithful \cite{big2}, \cite{kra}. Encouraged by these results, we ask ourselves if (new) link invariants can be obtained by similar methods.

Very briefly, the Lawrence representations of braid groups are constructed in two steps. Firstly, the braid group acts on a homotopy type called configuration space. Secondly, certain homology modules of configuration space are braid group modules.

In the case of links we expect the same two steps:
\begin{itemize}
\item From links to groups.
\item From groups to homology.
\end{itemize}
Our main result belongs to the first bullet. On the second bullet we have only some simple remarks.

Let $L\subset\rr^3$ be a link (not up to isotopy!) and fix a positive integer $k$. Consider the configuration spaces
\begin{align*}
C(L) &= \Big\{ X\subset\rr^3\backslash L\ \Big|\ |X|=k \Big\} \\
M(L) &= \Big\{ X\in C(L)\ \Big|\ \text{$X$ lies in a horizontal plane}\Big\}.
\end{align*}
It is trivial that up to diffeomorphism $C(L)$ depends only on the isotopy class of $L$. Therefore, a (twisted) homology module of $C(L)$ is a link invariant. But the fundamental group of $C(L)$ has no representations $U$ such that $H_*(C(L),U)$ is any interesting, at least no more than $\pi_1(\rr^3\backslash L)$ has.\footnote{This is because $\pi_1 C(L)$ is a semi-direct product $S_k\ltimes\pi_1(\rr^3\backslash L)^k$.} The group $\pi_1 M(L)$ has many more representations and we will henceforth concentrate on this group.

A particular case of our main result~\ref{st5} states:
\begin{quote} \em
If $L$ is a Heegaard link then $\pi_1 M(L)$ depends only on the isotopy class of $L$.
\end{quote}
(The full result considers the more general case of tangles.) See section~\ref{st35} for the definition of Heegaard links. Every link is isotopic to a Heegaard link.

A direct consequence of the above result is the construction of a link invariant which takes isomorphism classes of groups for values.

I don't know which other properties of $M(L)$ depend only on the isotopy class of $L$ ($L$ again Heegaard). Does $M(L)$ up to diffeomorphism? Does it up to homotopy equivalence?

The paper is built as follows. In section~\ref{st34} we review Lawrence's representations. Tangles and Heegaard tangles are introduced in section~\ref{st35}. The main result is formulated in section~\ref{st33} and proved in section~\ref{st36}. Section~\ref{st37} discusses the second bullet (from groups to homology) but it doesn't get very far.

\tableofcontents

\section{Lawrence representations} \label{st34}

We will review Lawrence's representations of braid groups \cite{law}.

The {\em braid group\,} $B_n$ is defined to be the fundamental group of
\[ BS_n=\big\{ X\subset\cc : |X|=n\big\}, \]
the space of sets (called {\em configurations}) of $n$ complex numbers.

Throughout this paper, we fix a natural number $k\geq 1$. Let $E$ denote the space of pairs $(X,Y)$ where $X\in BS_n$ and $Y\subset\cc\backslash X$ is a set of $k$ points which avoid $X$. The map
\begin{align*}
f\col E &\longrightarrow BS_n \\
(X,Y) &\longmapsto X
\end{align*}
is a fibre bundle (topologically locally trivial map). Let $F\subset E$ be the fibre of $f$ over a base-point in $BS_n$.

The fibre bundle $f$ admits a continuous section
\[ s\col BS_n\longrightarrow E \]
(section means $fs(X)=X$ for all $X$), for example
\[ s(X)=\big(X,\{a_X+1,a_X+2,\ldots,a_X+k\}\big) \]
where $ a_X=\max\big\{|x| : x\in X\big\}$. It gives rise to a splitting \[ t\col B_n\rightarrow\pi_1 E. \]
The theory of fibre bundles says that $\pi_1 E$ is now a semi-direct product $B_n\ltimes\pi_1 F$. In particular we have an action
\be B_n\longrightarrow \Aut(\pi_1 F) \label{st25} \ee
defined by $x(y)=(tx)y(tx)^{-1}$ ($x\in B_n$, $y\in\pi_1 F$).

Let $U$ be a module over any ring and let $r\col\pi_1 E\rightarrow \GL(U)$ be a linear representation. We then put
\[ V:=H_k(F,U). \]
It is known that $F$ is a $K(\pi,1)$, so that we also have $V=H_k(\pi_1 F,U)$. The $B_n$-action (\ref{st25}) on $\pi_1 F$ gives rise to a $B_n$-action
\[ B_n\longrightarrow\GL(V) \]
because homology is a functor. This is a general form of the Lawrence representation of the braid group \cite{law}.

\begin{figure}
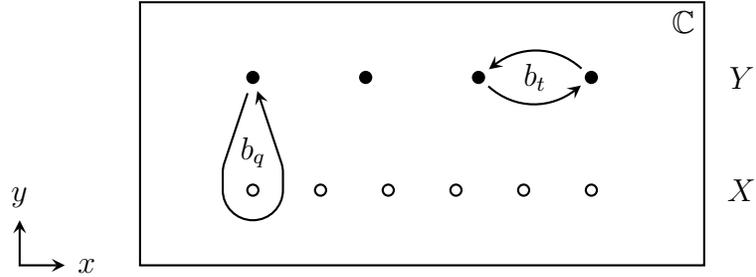

\caption{$b_q$ and $b_t$ \label{st26}}
\[
\pspicture(-1.5,0)(8.5,3.5)
\psset{unit=1mm}
\pspolygon(5,0)(80,0)(80,35)(5,35)
\psset{linearc=3.99}
\multido{\i=20+15}{4}{\pscircle*(\i,25){2.5pt}}
\multido{\i=20+9}{6}{\pscircle(\i,10){2.5pt}}
\psset{arrowsize=1pt 5, arrowlength=1, nodesep=5pt, arcangleA=-42, arcangleB=-38, ncurv=.8}
\pcarc{->}(50,25)(65,25) \pcarc{->}(65,25)(50,25)
\psline{->}(19.3,22.9)(16,13)(16,6)(24,6)(24,13)(20.6,23.2)
\rput(57.5,25){$b_t$}  \rput(20,15){$b_q$}
\uput[-135]{0}(80,35){$\cc$}
\rput(85,25){$Y$}   \rput(85,10){$X$}
\psset{linearc=0pt} \psline{<->}(-11,6)(-11,0)(-5,0)
\uput[90]{0}(-11,6){$y$} \uput[0]{0}(-5,0){$x$}
\endpspicture
\]
\end{figure}

The case where $U$ is $1$-dimensional is especially interesting. In that case, the space of representations $r$ is $2$-dimensional if $k\geq 2$; indeed, the abelianisation $(\pi_1 E)^{{\rm ab}}$ of $\pi_1 E$ is isomorphic to $\zz^2$. Figure~\ref{st26} shows two generators of $(\pi_1 E)^{{\rm ab}}$. White dots are elements of $X$ and shall be called {\em punctures}. Black dots are elements of $Y$. So $b_q\in\pi_1 E$ means that an element of $Y$ makes a full circle around an element of $X$ in counterclockwise direction, and $b_t\in\pi_1 E$ means that two elements of $Y$ interchange counterclockwise. Then $(\pi_1 E)^{{\rm ab}}$ is the free abelian group on the images of $b_q$ and $b_t$. The representation $r\col\pi_1 E\rightarrow\GL\big(1,\zz[q^{\pm1},t^{\pm1}]\big)$ is defined by $r(b_q)=q$ and $r(b_t)=t$. It can be shown that
\[ \dim V=\binom{n+k-2}{k} \]
so that the Lawrence representation can briefly be written
\[ \textstyle B_n\longrightarrow \GL(V)=\GL\Big(\binom{n+k-2}{k},\zz[q^{\pm1},t^{\pm1}]\Big). \]
The Jones polynomial has been related to these representations in \cite{law} and \cite{big3}. For $k=1$ this representation is the well-known Burau representation discovered in 1936. The representation for $k=2$ was shown to be faithful in \cite{big2} and \cite{kra}. 

In the following we will try to apply similar methods to obtain knot and link invariants.

\section{Tangles} \label{st35}

We define tangles and Heegaard tangles.

A {\em tangle\,} of type $[a,b]$ ($a,b\in\rr$, $a<b$) is a smooth compact 1-manifold $T\subset\cc\times[a,b]$ (with coordinates $(x+iy,z)$) with $\partial T=T\cap\big(\cc\times\{a,b\}\big)$ and such that $T$ is not tangent to $\cc\times\{a,b\}$. A {\em link\,} is a tangle with empty boundary.

Two tangles of types $[a,b]$ and $[c,d]$ are {\em isotopic\,} if one is taken to the other by a diffeomorphism $f\col\cc\times[a,b]\rightarrow\cc\times[c,d]$ with $f(x,a)=(x,c)$ and $f(x,b)=(x,d)$, for all $x\in\cc$. The isotopy class of a tangle $T$ is written $[T]$.

It may happen that the union $T_1\cup T_2$ of a tangle $T_1$ of type $[a,b]$ and a tangle $T_2$ of type $[b,c]$ is again a tangle (of type $[a,c]$); if this happens we call $T_1T_2:=T_1\cup T_2$ the {\em product\,} of $T_1$, $T_2$.

The isotopy class $[T_1T_2]$ depends only on $[T_1]$, $[T_2]$ and we thus obtain the {\em multiplication\,} of isotopy classes of tangles.

\begin{defi} Let
\begin{align*} 
p_3\col\cc\times\rr &\longrightarrow\rr \\
(x+iy,z) &\longmapsto z
\end{align*}
denote the projection on the third real coordinate. A plane $p_3^{-1}(d)$ ($d\in\rr$) is called {\em horizontal}. A tangle $T$ of type $[a,b]$ is said to be {\em Heegaard\,} if a number $c\in[a,b]$ is distinguished such that any local maximum (or {\em cap}) $x$ of
\[ p_3|_T\col T\longrightarrow \rr \]
satisfies $p_3(x)>c$, and every local minimum (or {\em cup}) satisfies $p_3(x)<c$, and $T$ is not tangent to $p_3^{-1}(c)$. The horizontal plane $p_3^{-1}(c)$ is called the {\em Heegaard plane\,} and it separates the caps from the cups.
\end{defi}

Two Heegaard tangles are said to be {\em Heegaard isotopic\,} if one is taken to the other by a diffeomorphism $f\col\cc\times[a,b]\rightarrow\cc\times[c,d]$ as for the usual isotopy provided $f$ takes the one Heegaard plane to the other.

If two Heegaard tangles are Heegaard isotopic then they are isotopic, but not conversely. It is known \cite{bir} that every tangle is isotopic to a Heegaard tangle.

The concept of Heegaard tangles is closely related to the {\em plat closure\,} about which we shall be rather brief. A detailed discussion can be found in \cite{bir}. We have a commuting diagram
\be 
\begin{aligned}
\begin{psmatrix}[colsep=4em]
\displaystyle\smash{\coprod_{n\geq 0}B_{2n}} && \Big\{ \text{Links} \Big\}\Big/ \text{isotopy} \\
&\makebox[0mm]{\mbox{$ \Big\{ \text{Heegaard links} \Big\} \Big/ \text{Heegaard isotopy}$}}
\psset{arrowsize=2pt 5, arrowlength=1}
\ncline[nodesep=1ex]{->}{1,1}{1,3}
\bput{0}{\text{\small plat closure}}
\ncline[nodesep=3ex]{->}{1,1}{2,2}
\ncline[nodesepA=4.5ex, nodesepB=1ex]{->}{2,2}{1,3}
\end{psmatrix}
\end{aligned} \label{st32}
\ee
and the map in the top row, the plat closure, adds $n$ caps and $n$ cups to a braid on $2n$ strings. Loosely speaking, the plat closure of a braid is in Heegaard position in a natural way. The elements of the set at the bottom of (\ref{st32}) can be viewed as certain double cosets of braids as is done in \cite{bir}.

The term plat closure is well-known, but we prefer the language of Heegaard tangles because they are not isotopy classes, contrary to plat closures.

\section{From tangles to groups} \label{st33}

\begin{defi} Let $T$ be a tangle of type $[a,b]$. We define
\[ M(T)=\Big\{ X\subset\big(\cc\times[a,b]\big)\backslash T\ \Big|\ x,y\in X\Rightarrow p_3(x)=p_3(y) \Big\} \]
which we call the {\em configuration space\,} of $T$. Note that every $X\in M(T)$ is required to lie in a horizontal plane.
\end{defi}

For example, if $k=1$ then $M(T)=\big(\cc\times[a,b]\big)\backslash T$.

\begin{defi} \label{st1} A Heegaard tangle is {\em saturated\,} if each of its components contains a cap or a cup.
\end{defi}

It is known \cite{bir} that every tangle is isotopic to a Heegaard tangle, and therefore clearly also to a saturated Heegaard tangle. Every Heegaard link is saturated. 

The following is our main result.

\begin{theo} \label{st5} Let $T_1,T_2$ be saturated Heegaard tangles. If $T_1,T_2$ are (`non-Heegaard') isotopic then
\[ \pi_1 M(T_1)\cong \pi_1 M(T_2). \]
\end{theo}
 
The point of the theorem is that $T_1$, $T_2$ are not assumed to be Heegaard isotopic  but just isotopic, which is a weaker assumption. It is trivial that $\pi_1 M(T_1)$ and $\pi_1 M(T_2)$ are isomorphic if $T_1,T_2$ are Heegaard isotopic.

A {\em tangle invariant\,} is just a map from the set of isotopy classes of tangles to any set. Theorem~\ref{st5} suggests a tangle invariant as follows.

\begin{defi} Let $T$ be a tangle. We define a group $F(T)$ as follows. First, choose a saturated Heegaard tangle $U$ isotopic to $T$. We put
\[ F(T)=\pi_1 M(U). \]
\end{defi}

By \ref{st5}, $F(T)$ is independent of the choice of $U$. It is immediate that $F(T)$ depends only on the isotopy class of $T$, so $T\mapsto F(T)$ is a tangle invariant.

I don't know which other properties of $M(T)$ ($T$ a saturated Heegaard tangle) depend only on the isotopy class of $T$. Is $M(T)$ up to diffeomorphism a tangle invariant? Is its homotopy type?

\section{Proof} \label{st36}

After some preparation, we will prove our main result~\ref{st5}.

In this section, all Heegaard tangles will be of type $[-1,1]$ and with Heegaard plane $H:=p_3^{-1}(0)$, unless stated otherwise.

\begin{figure}
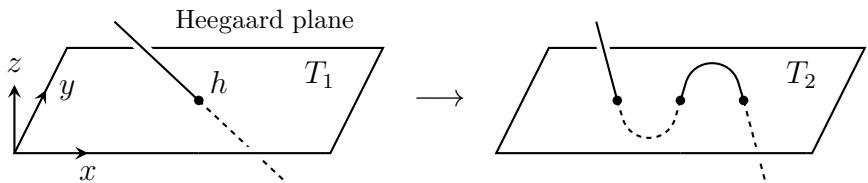

\caption{Elementary stabilisation \label{st27}}
\[
\pspicture[.5](-2.45,-1.1)(2.45,1.2)
\psset{unit=.7mm}
\pspolygon(25,-10)(35,10)(-25,10)(-35,-10)
\psset{arrowsize=1pt 5, arrowlength=1}
\psline[linestyle=none]{->}(-35,-10)(-21,-10)
\uput[-90]{0}(-21,-10){$x$}
\psline[linestyle=none]{->}(-35,-10)(-29,2)
\uput[0]{0}(-29,2){$y$}
\psline{->}(-35,-10)(-35,3)
\uput[90]{0}(-35,3){$z$}
\psline[border=2pt](-16,15)(0,0)
\psset{linestyle=dashed, dash=2pt 2pt}
\psline(16,-15)(0,0)
\psdots(0,0)       \rput(23,5){$T_1$}
\uput[40]{0}{$h$}  \psset{linestyle=solid}
\psline[border=2pt](0,-10)(24,-10)
\psset{linestyle=none}
\pcline(-20,10)(45,10) \Aput{{\footnotesize Heegaard plane}}
\endpspicture
\raisebox{-.7ex}{$\quad\longrightarrow\quad$}
\pspicture[.5](-2.45,-1.1)(2.45,1.2)
\psset{unit=.7mm}
\pspolygon(25,-10)(35,10)(-25,10)(-35,-10)
\psset{linearc=5.2}
\psline[border=2pt](-16,15)(-12,0)
\psline(0,0)(6,20)(12,0)
\psset{linestyle=dashed, dash=2pt 2pt}
\psline(-12,0)(-6,-20)(0,0)  \psline(12,0)(16,-15)
\psdots(-12,0)(0,0)(12,0)    \rput(23,5){$T_2$}
\psset{linestyle=solid}
\psline[border=2pt](0,-10)(24,-10)
\endpspicture
\]
\end{figure}

\kopje{Stabilisation}
Let $T_1,T_2$ be Heegaard tangles. We say that $T_2$ is obtained from $T_1$ by an {\em elementary stabilisation\,} if $T_1,T_2$ only differ close to some intersection point $h\in T_1\cap H$ and $T_2$ has three such intersection points close to $h$ rather than one. See figure~\ref{st27}.

Note that in the foregoing, $T_2$ is determined up to Heegaard isotopy by $T_1$ and $h$ (and isotopy which is trivial far away from $h$).

We will use the following result by Birman \cite{bir}. In fact she considered only links, but her proof also works for tangles.

\begin{theo} The obvious map
\[ \Big\{ \text{\rm Heegaard links} \Big\} \Big/
\substack{ \text{\rm Heegaard isotopy} \\[.2ex] \text{\rm  and stabilisation}} \longrightarrow
\Big\{ \text{\rm Links} \Big\}\Big/\text{\rm isotopy} \]
is bijective. (The set on the left is by definition the quotient of the set of Heegaard links by the equivalence relation $\sim$ generated by $T_1\sim T_2$ whenever $T_1,T_2$ either are Heegaard isotopic or differ by an elementary stabilisation.)\blokk
\end{theo}

\kopje{An application of the Seifert-Van Kampen theorem}
Let $T$ be a Heegaard tangle. Let $Z$ denote the map
\[  \begin{aligned} Z\col M(T) &\longrightarrow \rr \\ X &\longmapsto p_3(x)\text{ for one (hence all) $x\in X$.}
\end{aligned} \]
We will use the following notation.
\begin{align*}
M &= M(T) &G&=G(T)=\pi_1 M(T) \\
M_0 &= M_0(T) = Z^{-1}(0) &G_0&=G_0(T)=\pi_1 M_0(T) \\
M_{+} &= M_{+}(T) = Z^{-1}\big([0,1]\big) &G_{+}&=G_{+}(T)=\pi_1 M_{+}(T) \\
M_{-} &= M_{-}(T) = Z^{-1}\big([-1,0]\big) &G_{-}&=G_{-}(T)=\pi_1 M_{-}(T)
\end{align*}

An immediate application of the Seifert-Van Kampen theorem shows that the diagram
\be 
\begin{aligned} 
\begin{psmatrix}[colsep=3ex, rowsep=2ex]
& G_{+} \\
G_0 && G \\
& G_{-} 
\psset{arrowsize=2pt 5, arrowlength=1, nodesep=.6ex, linestyle=none}
\ncline{-}{2,1}{1,2} \lput{:U}{\longsur}
\ncline{-}{2,1}{3,2} \lput{:U}{\longsur}
\ncline{-}{1,2}{2,3} \lput{:U}{\longsur}
\ncline{-}{3,2}{2,3} \lput{:U}{\longsur}
\end{psmatrix} 
\end{aligned}
\label{st2} \ee
(with obvious arrows) is a push-out diagram. There is another way of saying the same thing, because all maps in (\ref{st2}) are surjective: writing
\begin{align}
K_{+}&=K_{+}(T)=\ker(G_0\rightarrow G_{+}) \notag \\
K_{-}&=K_{-}(T)=\ker(G_0\rightarrow G_{-}) \notag \\
K &= K(T) =\ker(G_0\rightarrow G) \label{st39}
\end{align}
we have
\be K=\<K_{+},K_{-}\>. \label{st11} \ee
Each of the three straight lines of two arrows in
\[ \begin{psmatrix}[colsep=3ex, rowsep=3ex]
& K_{-} && G_{+} \\
K && G_0 && G \\
& K_{+} && G_{-} 
\psset{arrowsize=2pt 5, arrowlength=1, nodesep=.6ex, linestyle=none}
\ncline{-}{2,1}{1,2} \lput{:U}{\longinj}
\ncline{-}{2,1}{3,2} \lput{:U}{\longinj}
\ncline{-}{1,2}{2,3} \lput{:U}{\longinj}
\ncline{-}{3,2}{2,3} \lput{:U}{\longinj}
\ncline{-}{2,1}{2,3} \lput{:U}{\hugeinj}
\ncline{-}{2,3}{1,4} \lput{:U}{\longsur}
\ncline{-}{2,3}{3,4} \lput{:U}{\longsur}
\ncline{-}{1,4}{2,5} \lput{:U}{\longsur}
\ncline{-}{3,4}{2,5} \lput{:U}{\longsur}
\ncline{-}{2,3}{2,5} \lput{:U}{\hugesur}
\end{psmatrix}
 \]
is exact.

\kopje{Some group presentations}
The following lemma gives generators for $K_{+}(T)$.

\begin{figure}
\caption{$T\cap\big(\cc\times[0,1]\big)$ \label{st29}}
\[
\pspicture(-5,-2)(5,2.5)
\psset{unit=1.25mm}
\pspolygon(20,-15)(40,15)(-20,15)(-40,-15)
\psellipse(-20,-6)(10,5)  \psellipse(0,-2)(10,5)
\psellipse(18,6)(10,5)    \psellipse(-10,8)(10,5)
\psset{nodesep=0pt, arcangleA=90, arcangleB=90, ncurv=1, showpoints=true}
\pcarc{-}(-25,-6)(-15,-6)  \pcarc{-}(13,6)(23,6)
\psset{arcangleA=45, arcangleB=135, ncurv=1.4}
\pcarc{-}(-2.2,-4.2)(2.2,0.2)
\psset{arcangleA=120, arcangleB=60, ncurv=1.2}
\pcarc{-}(-13.5,10)(-6.5,6)
\psdots(-25,-6)(-15,-6)(13,6)(23,6)(-2.2,-4.2)(2.2,0.2)(-13.5,10)(-6.5,6)
\psset{border=2pt}
\psline(-27,0)(-27,15)  \psline(16,-2)(16,13)
\psline(30,5)(30,20)    \psline(20,-10)(20,5)
\psline(20,-10)(20,5)   \psline(-18,2)(-18,17)
\psline(5,-10)(5,5)
\uput[0]{0}(30,0){$H$}  \rput(-20,-7){$D_i$}
\multido{\i=2+2}{4}{\pscircle*(\i,12){2.5pt}}
\psset{arrowsize=1pt 5, arrowlength=1, border=0pt, showpoints=false}
\psline[linestyle=none]{->}(-40,-15)(-32,-15)
\uput[-90]{0}(-32,-15){$x$}
\psline[linestyle=none]{->}(-40,-15)(-36,-9)
\uput[0]{0}(-36,-9){$y$}
\psline{->}(-40,-15)(-40,-8)
\uput[90]{0}(-40,-8){$z$}
\endpspicture
\]
\end{figure}

\begin{lemm} \label{st3} Let $T$ be a Heegaard tangle. Let $D_1,\ldots,D_\ell$ be disjoint closed disks in the Heegaard plane $H$ and write $\overline{D_i}:=p_3 D_i$. Suppose that
\begin{align*}
\big(\partial \overline{D_i} \times[0,1]\big) \cap T &= \varnothing\qquad\text{for all $i$}
\end{align*}
and that $\big( \overline{D_i}\times[0,1]\big) \cap T$ is an interval whose boundary lies in $H$, for all~$i$. See figure~\ref{st29}. Loosely speaking every cap of $T$ lives in another $\overline{D_i}\times[0,1]$. Recall that the Heegaard plane is $H=p_3^{-1}(0)$ so we are only looking at the part above $H$.

Fix a set $X_0\subset H\backslash(D_1\cup\cdots \cup D_\ell)$ of $k-1$ elements (the heavy dots in figure~\ref{st29}). Consider the conjugacy class $Y_i\subset G_0(T)$ of those elements given by a closed path in $M_0(T)$ homotopic to the map
\[ \left\{ \begin{aligned}
\partial D_i &\longrightarrow M_0(T) \\
x &\longmapsto \{x\}\cup X_0. 
\end{aligned} \right. \]
Then $K_{+}(T)$ is generated by $Y_1\cup\cdots\cup Y_\ell$.
\end{lemm}

\proof Left to the reader.\blokk

\begin{rema} One can get around the need of Lemma~\ref{st3} if one replaces $\pi_1 M(T)$ in the main result \ref{st5} by the group implied in Lemma~\ref{st3}, that is,
\[ G_0(T)/\<Y_1,\ldots,Y_\ell,Z_1,\ldots,Z_m\> \]
where the $Y_i$ are as in \ref{st3} and $Z_i$ likewise with cups instead of caps. The price one pays is that one should show that this group is well defined, that is, does not depend on the choice of the disks $D_i$ in \ref{st3}. Modification of the proof of \ref{st5} is not necessary.
\end{rema}

\begin{figure}
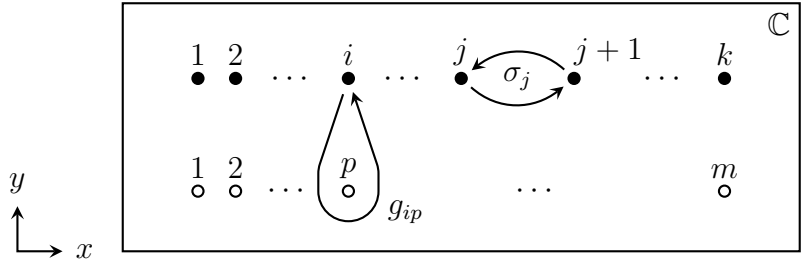

\caption{Generators for the braid group of the punctured disk \label{st28}}
\[
\pspicture(-2.4,0.2)(8,3.5)
\psset{unit=1mm}
\pspolygon(-10,2)(80,2)(80,35)(-10,35)
\psset{linearc=3.99}
\pscircle*(0,25){2.5pt}   \uput[90]{0}(0,25){$1$}
\pscircle*(5,25){2.5pt}   \uput[90]{0}(5,25){$2$}
\pscircle*(20,25){2.5pt}  \uput[90]{0}(20,25){$i$}
\pscircle*(35,25){2.5pt}  \uput[90]{0}(35,25){$j$}
\pscircle*(50,25){2.5pt}  \uput[90]{0}(50,25){\makebox[0mm][l]{$j+1$}}
\pscircle*(70,25){2.5pt}  \uput[90]{0}(70,25){$k$}
\pscircle(0,10){2.5pt}    \uput[90]{0}(0,10){$1$}
\pscircle(5,10){2.5pt}    \uput[90]{0}(5,10){$2$}
\pscircle(20,10){2.5pt}   \uput[90]{0}(20,10){$p$}
\pscircle(70,10){2.5pt}   \uput[90]{0}(70,10){$m$}
\psset{arrowsize=1pt 5, arrowlength=1, nodesep=5pt, arcangleA=-42, arcangleB=-38, ncurv=.8}
\pcarc{->}(35,25)(50,25) \pcarc{->}(50,25)(35,25)
\psline{->}(19.3,22.9)(16,13)(16,6)(24,6)(24,13)(20.6,23.2)
\rput(12.5,25){$\cdots$}  \rput(27.5,25){$\cdots$}
\rput(62,25){$\cdots$}    \rput(12,10){$\cdots$}
\rput(45,10){$\cdots$}    \uput[-135]{0}(80,35){$\cc$}
\rput(42.5,25){$\sigma_j$}
\uput[-45]{0}(24,10){$g_{ip}$}
\psset{linearc=0pt} \psline{<->}(-24,8)(-24,2)(-18,2)
\uput[90]{0}(-24,8){$y$} \uput[0]{0}(-18,2){$x$}
\endpspicture
\]
\end{figure}

\begin{prop} \label{st12} The $k$-string braid group 
\[ \pi_1\Big\{ X\subset\cc \ \Big|\ |X|=k,\ X\cap\{1,\ldots,m\}=\varnothing \Big\} \]
of the $m$ times punctured disk with base-point $\{1+i,2+i,\ldots,k+i\}$ ($i=\sqrt{-1}$) is presented by generators
\be
\begin{aligned}
\sigma_i &\qquad&& 1\leq i< k \\
g_{ip} &&& 1\leq i\leq k,\ 1\leq p\leq m
\end{aligned}
\label{st13} \ee
(see figure~\ref{st28}) and relations
\begin{align}
\label{st6} \sigma_i \sigma_j &= \sigma_j \sigma_i && |i-j|>1 && \\
 \sigma_i \sigma_j \sigma_i &= \sigma_j \sigma_i \sigma_j && |i-j|=1 \label{st7} \\
\label{st8} [g_{ip},g_{jq}] &=1 &&i<j,\ p<q \\
\label{st9} \hspace*{2em} \sigma_i\, g_{i+1,p}\, \sigma_i^{-1}\, g_{ip}^{-1} &= 1 \\
\label{st10} [g_{ip},\sigma_i\, g_{ip}\, \sigma_i] &= 1 \\
\label{st31} [g_{ip},\sigma_j] &=1 &&i\neq\{j,j+1\}
\end{align}
where $[a,b]=aba^{-1}b^{-1}$.
\end{prop}

\proof Presentations for this group can be found in Theorem~2 or Theorem~3 of \cite{lam} or Theorem~5.1 of \cite{bel}. The generators $\sigma_j$, $a_{ij}$ of Theorem~2 in \cite{lam} are our $\sigma_j$, $g_{j,m+1-i}^{-1}$. It is left to the reader to check that this identification respects the group presentations.\blok

For any tangle $T$, the group $G_0(T)$ is the $k$-string braid group of a punctured disk $H\backslash T$.

\begin{figure}
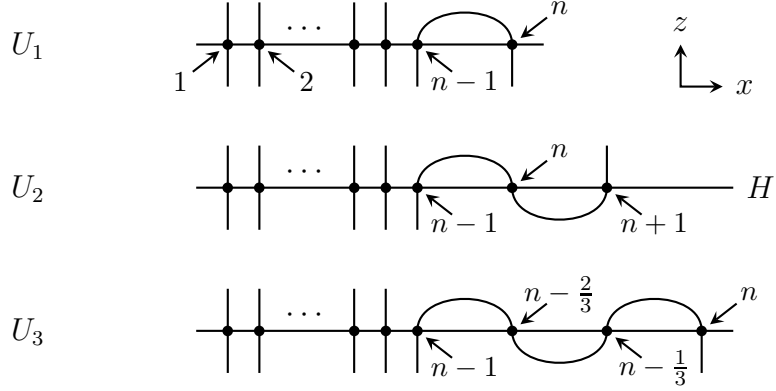

\caption{Three fundamental tangles \label{st30}}
\begin{align*}
\raisebox{-.7ex}{$U_1$}
&&&
\pspicture[.5](-3,-.7)(4.4,.7)  \SpecialCoor
\psset{yunit=.7mm, xunit=.7mm}
\psline(0,-8)(0,0)   \psline(18,-8)(18,0)
\psset{arcangleA=90, arcangleB=90, ncurv=.9}
\pcarc{-}(0,0)(18,0)   \psline(-12,-8)(-12,8)
\psline(-6,-8)(-6,8)   \psline(-42,0)(24,0)
\rput(-21,3){$\cdots$}
\psline(-36,-8)(-36,8) \psline(-30,-8)(-30,8)
\psset{dotsize=4pt}
\psdots(0,0)(18,0)(-12,0)(-6,0)(-36,0)(-30,0)
\psset{arrowsize=1pt 5, arrowlength=1}
\rput(9,-7){\rnode{v1}{{\small $n-1$}}}
\pcline[nodesepA=6pt, nodesepB=4pt]{->}(v1)(0,0)
\rput(-45,-7){\rnode{w1}{{\small $1$}}}
\pcline[nodesepA=6pt, nodesepB=4pt]{->}(w1)(-36,0)
\rput(-21,-7){\rnode{w2}{{\small $2$}}}
\pcline[nodesepA=6pt, nodesepB=4pt]{->}(w2)(-30,0)
\rput(27,7){\rnode{v2}{{\small $n$}}}
\pcline[nodesepA=6pt, nodesepB=4pt]{->}(v2)(18,0)
\psline{<->}(50,0)(50,-8)(58,-8)
\uput[90]{0}(50,0){$z$} \uput[0]{0}(58,-8){$x$}
\endpspicture \\[2ex]
\raisebox{-.7ex}{$U_2$}
&&&
\pspicture[.5](-3,-.7)(4.4,.7)  \SpecialCoor
\psset{yunit=.7mm, xunit=.7mm}
\psline(0,-8)(0,0)   \psline(36,8)(36,0)
\psset{arcangleA=90, arcangleB=90, ncurv=.9}
\pcarc{-}(0,0)(18,0)    \pcarc{-}(36,0)(18,0)   
\psline(-12,-8)(-12,8)  \psline(-6,-8)(-6,8)   
\psline(-42,0)(60,0)    \rput(-21,3){$\cdots$}
\psline(-36,-8)(-36,8)  \psline(-30,-8)(-30,8)
\uput[0]{0}(60,0){$H$}
\psset{dotsize=4pt}
\psdots(0,0)(18,0)(-12,0)(-6,0)(-36,0)(-30,0)(36,0)
\psset{arrowsize=1pt 5, arrowlength=1}
\rput(9,-7){\rnode{v1}{{\small $n-1$}}}
\pcline[nodesepA=6pt, nodesepB=4pt]{->}(v1)(0,0)
\rput(27,7){\rnode{v2}{{\small $n$}}}
\pcline[nodesepA=6pt, nodesepB=4pt]{->}(v2)(18,0)
\rput(45,-7){\rnode{v3}{{\small $n+1$}}}
\pcline[nodesepA=6pt, nodesepB=4pt]{->}(v3)(36,0)
\endpspicture \\[2ex]
\raisebox{-.7ex}{$U_3$}
&&&
\pspicture[.5](-3,-.7)(4.2,.7)  \SpecialCoor
\psset{yunit=.7mm, xunit=.7mm}
\psline(0,-8)(0,0)   \psline(54,-8)(54,0)
\psset{arcangleA=90, arcangleB=90, ncurv=.9}
\pcarc{-}(0,0)(18,0)    \pcarc{-}(36,0)(18,0)   
\pcarc{-}(36,0)(54,0)   
\psline(-12,-8)(-12,8)  \psline(-6,-8)(-6,8)   
\psline(-42,0)(60,0)    \rput(-21,3){$\cdots$}
\psline(-36,-8)(-36,8)  \psline(-30,-8)(-30,8)
\psset{dotsize=4pt}
\psdots(0,0)(18,0)(-12,0)(-6,0)(-36,0)(-30,0)(36,0)(54,0)
\psset{arrowsize=1pt 5, arrowlength=1}
\rput(9,-7){\rnode{v1}{{\small $n-1$}}}
\pcline[nodesepA=6pt, nodesepB=4pt]{->}(v1)(0,0)
\rput(27,7){\rnode{v2}{{\small $n-\frac{2}{3}$}}}
\pcline[nodesepA=6pt, nodesepB=4pt]{->}(v2)(18,0)
\rput(45,-7){\rnode{v3}{{\small $n-\frac{1}{3}$}}}
\pcline[nodesepA=6pt, nodesepB=4pt]{->}(v3)(36,0)
\rput(63,7){\rnode{v4}{{\small $n$}}}
\pcline[nodesepA=6pt, nodesepB=4pt]{->}(v4)(54,0)
\endpspicture 
\end{align*}
\end{figure}

We consider three Heegaard tangles $U_1$, $U_2$, $U_3$ defined by figure~\ref{st30}. We have
\begin{align*}
U_1\cap H &= \{1,\ldots,n\}\times\{0\} \\
U_2\cap H &= \{1,\ldots,n+1\}\times\{0\} \\
U_3\cap H &= \textstyle \{1,\ldots,n-1,n-\frac{2}{3},n-\frac{1}{3},n\}\times\{0\}. 
\end{align*}
In the next lemma, we will prove that $G(U_1)$ and $G(U_2)$ are isomorphic in a precise sense. Of course, $G_0(U_1)$ is just the braid group of the $n$ times punctured disk. By definition, $G(U_1)=G_0(U_1)/K(U_1)$. Combining (\ref{st11}) and Lemma~\ref{st12} then shows that $G(U_1)$ is presented by generators (\ref{st13}) and relations (\ref{st6})--(\ref{st31}) (with $m=n$) as well as
\be g_{i,n-1}\,g_{in}=1 \qquad \text{for all $i\in\{1,\ldots,k\}$}. \label{st16} \ee
Similarly, $G(U_2)$ is presented by generators (\ref{st13}) and relations (\ref{st6})--(\ref{st31}) (with $m=n+1$) and
\begin{align}
 g_{i,n-1}\,g_{in} &=1 \qquad \text{for all $i\in\{1,\ldots,k\}$} \label{st19} \\
 g_{in}\,g_{i,n+1} &=1 \qquad \text{for all $i\in\{1,\ldots,k\}$}. \label{st20}
\end{align}

\begin{lemm} \label{st24} There is a (unique) isomorphism $f\col G(U_2)\rightarrow G(U_1)$ such that
\begin{align*}
f(\sigma_i) &= \sigma_i \\
f(g_{ip}) &= g_{ip} \hspace*{3em} (p\neq n+1) \\
f(g_{i,n+1}) &= g_{i,n-1}.
\end{align*}
\end{lemm}

\proof We need to prove that the substitution
\be g_{i,n+1}\longmapsto g_{i,n-1} \label{st15} \ee
takes any relation for $G(U_2)$ to one of the relations for $G(U_1)$ or a consequence of them. (It is clear that all relations of $G(U_1)$ are obtained this way.)

First consider (\ref{st8}) with $q=n+1$. 

Suppose $p<n$. Then we have the following computation in $G(U_1)$:
\be [g_{ip},g_{j,n-1}]
\stackrel{(\ref{st19})}{=} [g_{ip},g_{j,n}^{-1}]
\stackrel{(\ref{st8})}{=} 1. \label{st38} \ee
The substitution (\ref{st15}) takes (\ref{st8}) to $[g_{ip},g_{j,n-1}]=1$ which is true in $G(U_1)$ by (\ref{st38}). 

Suppose $p=n+1$. Then the following holds in $G(U_1)$:
\be [g_{in},g_{j,n-1}]
\stackrel{(\ref{st16})}{=}  [g_{i,n-1}^{-1},g_{jn}^{-1}]
\stackrel{(\ref{st8})}{=}  1  \label{st18} \ee
But the substitution (\ref{st15}) takes (\ref{st8}) to $[g_{in},g_{j,n-1}]=1$ which is true by (\ref{st18}).

Our substitution takes (\ref{st9}) with $p=n+1$ to $\sigma_i\,g_{i+1,n-1}\,\sigma_{i}^{-1}\,g_{i,n-1}^{-1}=1$ which is true in $G(U_1)$ by (\ref{st9}). Likewise,
the substitution takes (\ref{st10}) with $p=n+1$ to
$[g_{i,n-1},\sigma_i\, g_{i,n-1}\, \sigma_i] = 1$ which is true in $G(U_1)$ by (\ref{st10}). Also, the substitution takes (\ref{st31}) with $p=n+1$ to $[g_{i,n-1},\sigma_j]=1$ which is true in $G(U_1)$ by (\ref{st31}).

The relation (\ref{st19}) is just (\ref{st16}). Our substition takes (\ref{st20}) to a void statement.

All relations for $G(U_2)$ that we haven't mentioned so far don't involve $g_{i,n+1}$ and are clearly taken to a relation for $G(U_1)$.\blok

We have 
\be H\backslash U_3\subset H\backslash U_1. \label{st21} \ee
Indeed, the difference between these two sets is precisely $\{n-\frac{2}{3},n-\frac{1}{3}\}\times\{0\}$. The inclusion (\ref{st21}) yields an inclusion $M_0(U_3)\subset M_0(U_1)$ which in turn induces a surjective map of their fundamental groups
\be p\col G_0(U_3)\longsur G_0(U_1). \label{st22} \ee

\begin{lemm} \label{st23} We have $G(U_3)\cong G(U_1)$ and $pK(U_3)=K(U_1)$.\end{lemm}

\proof One proves that there exists an isomorphism $q\col G(U_3)\rightarrow G(U_1)$ by applying \ref{st24} twice, once on $(U_1,U_2)$ and once on $(U_2,U_3)$. Inspection of \ref{st24} and the presentations of $G(U_i)$ also shows that $q$ can be taken to be induced by $p$. By the definition of $K(U_i)$ (\ref{st39}) it follows that $pK(U_3)= K(U_1)$.\blok

\noi{\bf Proof of the main result \ref{st5}.} Let $T_1$, $T_3$ be two saturated Heegaard tangles of type $[-2,2]$ with Heegaard plane $H=p_3^{-1}(0)$, and suppose that $T_1$, $T_3$ differ by an elementary stabilisation. Our aim is to prove $G(T_1)\cong G(T_3)$.

After applying a Heegaard isotopy to $T_1$ and $T_3$ if necessary, and changing the sign of the $z$-coordinate in both of them if necessary, we may assume
\begin{align*}
T_1\cap[-1,1]=U_1, &&
T_3\cap[-1,1]=U_3, &&
T_1\backslash U_1=T_3\backslash U_3
\end{align*}
where $U_1$ and $U_3$ are as in figure~\ref{st30}. (The cap in $T_1\cap[-1,1]$ exists because $T_1$ is saturated).

Note that
\begin{align*}
&&G_0(U_i) &=G_0(T_i) && (i=1,3) && 
\end{align*}
which therefore contains both $K(U_i)$ and $K(T_i)$. In Lemma~\ref{st23} we saw that the natural map $p$ in (\ref{st22}) takes $K(U_3)$ to $K(U_1)$. We will show that it also takes $K(T_3)$ to $K(T_1)$. 

Let $C$ denote the set of caps and cups of $T_1\backslash U_1=T_3\backslash U_3$. By \ref{st3} we have
\begin{align*}
&&K(U_i) &\subset K(T_i) && (i=1,3). && 
\end{align*}
Moreover one can associate, to each $c\in C$, two conjugacy classes $Y_i(c)\subset G_0(U_i)$ ($i=1,3$) such that
\begin{align*}
K(T_i) &= \big\< K(U_i),\{Y_i(c)\mid c\in C\}\big\>  \\
pY_3(c) &= Y_1(c)\qquad \text{for all $c\in C$.}
\end{align*}
We find
\begin{align*}
pK(T_3) &= p\big\< K(U_3),\{Y_3(c)\mid c\in C\}\big\>  \\
&= \big\< pK(U_3),\{pY_3(c)\mid c\in C\}\big\>  \\
&= \big\< K(U_1),\{Y_1(c)\mid c\in C\}\big\>  \\
&= K(T_1)
\end{align*}
as promised. Since the map $p$ is surjective we get
\[ G(T_3)=G_0(T_3)/K(T_3)\cong G_0(T_1)/K(T_1) = G(T_1). \]
This finishes the proof of our main result~\ref{st5}.\blok

\section{From groups to homology} \label{st37}

In section~\ref{st33} we defined a tangle invariant $F(T)$ which is a group. As we saw in the introduction, one hopes to turn this invariant into a more manageable invariant.

Suppose that we have, for each tangle $T$, a representation
\[ r(T)\col F(T)\longrightarrow \GL\big(U(T)\big) \]
defined over any ring. Suppose moreover that the pair $(F(T),r(T))$ is a tangle invariant up to isomorphism, in the sense that for any two isotopic tangles $T_1,T_2$ there exists a commutative diagram
\[ \begin{CD}
F(T_1) @>\makebox[3em][c]{{\scriptsize $r(T_1)$}}>> \GL\big(U(T_1)\big) \\
@VV\sim V @VV\sim V \\
F(T_2) @>\makebox[3em][c]{{\scriptsize $r(T_2)$}}>> \GL\big(U(T_2)\big)
\end{CD} \]
whose vertical arrows are isomorphisms, the right hand side one coming from an isomorphism $U(T_1)\rightarrow U(T_2)$. (Let us call a family $\{r(T)\}_T$ like this {\em good}.)  Then the isomorphism class of the homology module
\[ H_*\big(F(T),U\big) \]
is a tangle invariant.

For example, if $k=1$ and $U(T)$ is 1-dimensional then $H_1\big(F(T),U(T)\big)$ is the well-known coloured Alexander module.

For any tangle $T$ and any $k>0$ we have a continuous map
\begin{align*}
 M(T) &\longrightarrow BS_k \\
X &\longmapsto \big\{ x+iy\ \big|\ (x+iy,z)\in X\big\}
\end{align*}
(projection on the first two real coordinates) which induces a map to the braid group $G(T)\rightarrow B_k$. If also $T$ is a saturated Heegaard tangle then $F(T)=G(T)$ (see section~\ref{st33}) so that we have a map
\[ v(T)\col F(T)\longrightarrow B_k. \]
Every representation
\[ w\col B_k\longrightarrow\GL(U) \]
(not depending on any tangle) gives rise to a representation
\[ \begin{CD}
F(T) @>\makebox[2em][c]{{\scriptsize $v(T)$}}>> B_k @>\makebox[2em][c]{{\scriptsize $w$}}>> \GL(U)
\end{CD} \]
and therefore to a tangle invariant $H_*(F(T),U)$. This is one way to produce good families of representations $\{r(T)\}_T$ but certainly not the only way.

It would be interesting to compute any of the homology modules of $F(T)$ sketched in this section, or to know if they reveal information about links.


\end{document}